\documentclass[10pt]{article}
\usepackage{pb-diagram}
\usepackage{amsmath, amsfonts}
\usepackage{amssymb}
\usepackage{amscd}

\newtheorem{Df}{Definition}[section]
\newtheorem{Te}[Df]{Theorem}
\newtheorem{Po}[Df]{Proposition}
\newtheorem{Cr}[Df]{Corollary}
\newtheorem{Lm}[Df]{Lemma}
\newtheorem{Ca}[Df]{Claim}
\newtheorem{Cn}[Df]{Conjecture}
\newtheorem{Ex}[Df]{Example}
\newtheorem{Rm}[Df]{Remark}

\newcommand{\Bdf}{\begin{Df}}
\newcommand{\Edf}{\end{Df}}
\newcommand{\Bte}{\begin{Te}}
\newcommand{\Ete}{\end{Te}}
\newcommand{\Bpo}{\begin{Po}}
\newcommand{\Epo}{\end{Po}}
\newcommand{\Bcr}{\begin{Cr}}
\newcommand{\Ecr}{\end{Cr}}
\newcommand{\Blm}{\begin{Lm}}
\newcommand{\Elm}{\end{Lm}}
\newcommand{\Bca}{\begin{Ca}}
\newcommand{\Eca}{\end{Ca}}
\newcommand{\Bcn}{\begin{Cn}}
\newcommand{\Ecn}{\end{Cn}}
\newcommand{\Bex}{\begin{Ex}}
\newcommand{\Eex}{\end{Ex}}
\newcommand{\Brm}{\begin{Rm}}
\newcommand{\Erm}{\end{Rm}}

\begin{document}

\title{\bf{$N$-Cartan calculus}}
\author{Roland Berger}
\date{}

\maketitle

\begin{abstract}
We examine the $N$-Koszul calculus as defined by the author (Koszul calculus for $N$-homogeneous algebras, arXiv:1610.01035) in the example of the $N$-symmetric algebras. The case $N=2$ is fundamental in differential geometry since this case corresponds to the \'Elie Cartan calculus. Therefore, it is natural to conjecture the existence of an $N$-Cartan calculus on manifolds when $N>2$, which would provide a new prototype of noncommutative differential geometry. The goal of this paper is to express in terms of wedge product the defining formulas of the $N$-Koszul calculus for $N$-symmetric algebras. We hope that these formulas will be useful for performing an extension to manifolds. 
\end{abstract}

\section{The $N$-symmetric algebras}
In this paper, $k$ is a field, $V$ is a $k$-vector space of finite dimension $n$ with $n\geq 2$. For any $p\geq 1$, the $k$-linear map 
$Ant_p:V^{\otimes p}\rightarrow V^{\otimes p}$ is defined by
$$Ant_p(z)= \sum _{\sigma \in \Sigma_p} \mathrm{sgn}(\sigma)\, \sigma \cdot z,$$
where $\Sigma_p$ is the group of permutations of $\{1, \ldots , p\}$ and $\mathrm{sgn}$ is the signature. The left action $\sigma \cdot z$ of $\sigma$ on $z \in V^{\otimes p}$ is defined by
$$\sigma \cdot (x_1 \ldots x_p)= x_{\sigma^{-1} (1)} \ldots x_{\sigma^{-1} (p)},$$  
for any $x_1, \ldots, x_p$ in $V$. We omit the symbols $\otimes$ in the tensor products.  We set $Ant_0= Id_k$.

It is well-known that for any $p\geq 2$, we have 
\begin{equation} \label{imageant}
Ant_p(V^{\otimes p})=\bigcap_{i+2+j=p} V^{\otimes i}\otimes Ant_2(V^{\otimes 2}) \otimes V^{\otimes j}.
\end{equation}
Introducing the space $E$ generated by the elements $x\otimes x$ of $V^{\otimes 2}$ with $x$ in $V$, we have
\begin{equation} \label{kerant}
\ker (Ant_p)=\sum_{i+2+j=p} V^{\otimes i}\otimes E \otimes V^{\otimes j}.
\end{equation}
Therefore, the map $Ant_p$ induces a $k$-linear isomorphism
\begin{eqnarray} \label{isocan}
\bigwedge ^{p} V & \rightarrow & Ant_p(V^{\otimes p}) \\
e_{i_1} \wedge \ldots \wedge e_{i_p} & \mapsto & Ant_p(e_{i_1}  \ldots  e_{i_p}), \nonumber
\end{eqnarray}
where $1 \leq i_1 < \cdots < i_p \leq n$ and $(e_1, \ldots , e_n)$ is a basis of the space $V$. One has $Ant_p(V^{\otimes p})=0$ whenever $p>n$.

Throughout the paper, $N$ is an integer such that $2 \leq N \leq n$. We take $R=Ant_N(V^{\otimes N})$ as space of relations in the tensor algebra $T(V)$ of $V$. The corresponding $N$-homogeneous algebra
$$A=T(V)/(R)$$
is called \emph{the $N$-symmetric algebra of $V$} and is denoted by $A=S(V,N)$. For $N=2$, we recover the symmetric algebra $S(V,2)=S(V)$ of $V$.

Using the notation of~\cite{rb:Ncal}, we define the subspace
$$W_{p}=\bigcap_{i+N+j=p}V^{\otimes i}\otimes R\otimes V^{\otimes j}$$
of $V^{\otimes p}$ for any $p\geq 0$. According to (\ref{imageant}), we have
$$R=\bigcap_{i+2+j=N}V^{\otimes i}\otimes Ant_2(V^{\otimes 2}) \otimes V^{\otimes j},$$
and consequently
\begin{eqnarray} \label{spaceWp}
W_p =  & V^{\otimes p} & \mbox{if} \ 0 \leq p \leq N-1 \nonumber \\
W_p =  & Ant_p(V^{\otimes p}) & \mbox{if} \ N \leq p \leq n \\
W_p =  & 0 & \mbox{if} \ p > n. \nonumber 
\end{eqnarray}

Let us recall the following result~\cite{rb:nonquad, rbnm:kogo}. We do not know if the assumption on the characteristic of $k$ is necessary. 
\Bpo \label{Nko}
Assume that char(k) $=0$. The $N$-homogeneous algebra $A=T(V)/(R)$ is Koszul, of finite global dimension $d$, where 
\begin{enumerate}
\item $d = 2m$ \ \mbox{if} \ $n=Nm$,
  
\item $d = 2m+1$ \mbox{if} \ $n=Nm+r,\  0 < r < N$.
\end{enumerate}
If $N=2$, one has $d=n$ and $A$ is $n$-Calabi-Yau.

\noindent
If $N>2$, $A$ is AS-Gorenstein if and only if $n\equiv 1\ (N)$.

\noindent
If $N>2$, $A$ is $d$-Calabi-Yau if and only if $n\equiv 1\ (N)$ and $n$ is odd.
\Epo

Let us note that if $N>2$, $A$ contains a free algebra in two generators, so that the Gelfand-Kirillov dimension of $A$ is infinite and $A$ is not left (or right) noetherian.  

We will also need the map $\nu : \mathbb{N}\rightarrow \mathbb{N}$ defined by
\begin{eqnarray*}
\nu (p) &=& Np' \ \ \mathrm{if}\ p\ \mathrm{even},\ p=2p',\\
\nu (p) &=& Np'+1 \ \ \mathrm{if}\ p\ \mathrm{odd},\ p=2p'+1.
\end{eqnarray*}
Formulas (\ref{spaceWp}) are particularized in 
\begin{eqnarray} \label{spaceWnup}
W_{\nu(p)} =  & Ant_{\nu(p)}(V^{\otimes \nu(p)}) & \mbox{if} \ 0 \leq p \leq d \\
W_{\nu(p)} =  & 0 & \mbox{if} \ p > d. \nonumber 
\end{eqnarray}
Note that 
\begin{eqnarray*}
\nu (d) &=& n \ \ \mathrm{if}\ d\ \mathrm{even},\\
\nu (d) &=& Nm+1 \leq n \ \ \mathrm{if}\ d\ \mathrm{odd},\ d=2m+1.
\end{eqnarray*}

\setcounter{equation}{0}

\section{The Koszul cup product}

Let $P$ and $Q$ be $A$-bimodules where $A$ is an $N$-homogeneous algebra. For Koszul cochains $f:W_{\nu(p)}\rightarrow P$ and $g:W_{\nu(q)}\rightarrow Q$, we have defined in~\cite{rb:Ncal} the Koszul $(p+q)$-cochain $f\underset{K}{\smile} g : W_{\nu(p+q)} \rightarrow P\otimes_A Q$ by 
\\
1. if $p$ and $q$ are not both odd, so that $\nu(p+q)=\nu(p)+\nu(q)$, one has
\begin{equation*}
(f\underset{K}{\smile} g) (x_1 \ldots x_{\nu(p+q)}) = f(x_1 \ldots x_{\nu(p)})\otimes_A \, g(x_{\nu(p)+1} \ldots  x_{\nu(p)+\nu(q)}),
\end{equation*}
2. if $p$ and $q$ are both odd, so that $\nu(p+q)=\nu(p)+\nu(q)+N-2$, one has
\begin{eqnarray*}
    (f\underset{K}{\smile} g) (x_1 \ldots x_{\nu(p+q)})  =  -\sum_{0\leq i+j \leq N-2} x_1\ldots x_i f(x_{i+1} \ldots x_{i+\nu(p)})x_{i+\nu(p)+1}\ldots x_{\nu(p)+N-j-2} \\
  \otimes_A \ g(x_{\nu(p)+N-j-1} \ldots  x_{\nu(p)+\nu(q)+N-j-2})x_{\nu(p)+\nu(q)+N-j-1} \ldots  x_{\nu(p)+\nu(q)+N-2}.
\end{eqnarray*}

Let us rewrite these formulas for our example $A=S(V,N)$ in terms of wedge product. Suppose firstly that $\nu(p+q)=\nu(p)+\nu(q)$. Then the natural inclusion $W_{\nu(p+q)}\subseteq W_{\nu(p)}\otimes W_{\nu(q)}$ can be expressed by the standard equalities
\begin{eqnarray} \label{pqant1}
Ant_{\nu(p+q)}(x_1 \ldots x_{\nu(p+q)})= \sum _{\sigma \in Sh(\nu(p),\nu(q))} \mathrm{sgn}(\sigma) \,Ant_{\nu(p)}(x_{\sigma(1)} \ldots x_{\sigma(\nu(p))}) \nonumber \\
Ant_{\nu(q)}(x_{\sigma(\nu(p)+1)} \ldots x_{\sigma(\nu(p+q))}),   
\end{eqnarray}
where $Sh(\nu(p),\nu(q))$ denotes the subset of the $(\nu(p),\nu(q))$-shuffles in the group  $\Sigma_n$ of permutations of $\{1, \ldots , n\}$, and $x_1, \ldots ,x_{\nu(p+q)}$ are elements of the basis $(e_1, \ldots , e_n)$ of $V$ such that
$x_1< \ldots <x_{\nu(p+q)}$.

Using the isomorphisms (\ref{isocan}), we consider $f:\bigwedge ^{p} V \rightarrow P$ and $g:\bigwedge ^{q} V \rightarrow Q$. We arrive to 
\begin{eqnarray} \label{cup1}
(f\underset{K}{\smile} g) (x_1 \wedge \ldots \wedge x_{\nu(p+q)}) = (-1)^{pq} \sum _{\sigma \in Sh(\nu(p),\nu(q))} \mathrm{sgn}(\sigma) \,f(x_{\sigma(1)} \wedge \ldots \wedge x_{\sigma(\nu(p))})\nonumber \\
\otimes _A \, g(x_{\sigma(\nu(p)+1)} \wedge \ldots \wedge x_{\sigma(\nu(p+q))}).   
\end{eqnarray}

Suppose now that $\nu(p+q)\neq \nu(p)+\nu(q)$, meaning that $N>2$ and $p$, $q$ are both odd. For $0\leq i+j \leq N-2$, the natural inclusion
$$W_{\nu(p+q)}\subseteq V^{\otimes i} \otimes W_{\nu(p)}\otimes V^{\otimes N-2-i-j} \otimes W_{\nu(q)}\otimes V^{\otimes j}$$
can be expressed by the (less?) standard equalities
\begin{eqnarray} \label{pqant2}
  \lefteqn{Ant_{\nu(p+q)}(x_1 \ldots x_{\nu(p+q)}) =\sum _{\sigma \in Sh(i,\nu(p),\nu(q),j)} \mathrm{sgn}(\sigma) } \nonumber \\
& & x_{\sigma(1)} \ldots x_{\sigma(\nu(i))} Ant_{\nu(p)}(x_{\sigma(i+1)} \ldots x_{\sigma(i+\nu(p))})\, x_{\sigma(i+\nu(p)+1)} \ldots x_{\sigma(\nu(p)+N-j-2)} \\
& & Ant_{\nu(q)}(x_{\sigma(\nu(p)+N-j-1)} \ldots x_{\sigma(\nu(p)+\nu(q)+N-j-2)})\, x_{\sigma(\nu(p)+\nu(q)+N-j-1)} \ldots x_{\sigma(\nu(p+q))}, \nonumber
\end{eqnarray}
where $Sh(i,\nu(p),\nu(q),j)$ denotes the subset of $\Sigma_n$ formed of the $(\nu(p),\nu(q))$-shuffles acting on the index sets $\{i+1, \ldots ,i+\nu(p)\}$ and $\{\nu(p)+N-j-1, \ldots , \nu(p)+\nu(q)+N-j-2\}$. 

Using the isomorphisms (\ref{isocan}), we arrive to 
\begin{eqnarray} \label{cup2}
\lefteqn{(f\underset{K}{\smile} g) (x_1 \wedge \ldots \wedge x_{\nu(p+q)})= - \sum_{0\leq i+j \leq N-2} \ \ \sum _{\sigma \in Sh(i,\nu(p),\nu(q),j)} \mathrm{sgn}(\sigma)  } \nonumber \\
& & x_{\sigma(1)} \ldots x_{\sigma(\nu(i))} f(x_{\sigma(i+1)}\wedge \ldots \wedge x_{\sigma(i+\nu(p))})\, x_{\sigma(i+\nu(p)+1)} \ldots x_{\sigma(\nu(p)+N-j-2)} \\
& & \otimes_A g(x_{\sigma(\nu(p)+N-j-1)}\wedge \ldots \wedge x_{\sigma(\nu(p)+\nu(q)+N-j-2)})\, x_{\sigma(\nu(p)+\nu(q)+N-j-1)} \ldots x_{\sigma(\nu(p+q))}, \nonumber  
\end{eqnarray}
where $x_1< \ldots <x_{\nu(p+q)}$ in $(e_1, \ldots , e_n)$.
\\

If $N=2$ and $P=Q=A$, the cup product $f \underset{K}{\smile} g$ coincides with the Cartan wedge product $f\wedge g$, so that it is associative and graded commutative in this case. It would be interesting to examine whether associativity holds if $N>2$. The same for graded commutativity, but note that $A$ is not commutative when $N>2$.
\\ 

\textbf{Question 1.} As in the quadratic case $N=2$, is it possible to extend Formulas (\ref{cup1}) and (\ref{cup2}) to smooth functions -- when $k=\mathbb{R}$ and specializing cochains on $P$ and $Q$ to sections of vector bundles -- in such a way the extended formulas are compatible to chart changes? Shortly, is it possible to define an $N$-Koszul cup product on manifolds?

\setcounter{equation}{0}

\section{The Koszul cap products}

\subsection{Definition and first properties}

Let $A$ be an $N$-homogeneous algebra, $M$ and $P$ be $A$-bimodules. For any Koszul $p$-cochain $f:W_{\nu(p)}\rightarrow P$ and any Koszul $q$-chain $z=m \otimes x_1 \ldots x_{\nu(q)}$ 
in $M\otimes W_{\nu(q)}$, we have defined in~\cite{rb:Ncal} the Koszul $(q-p)$-chains $f \underset{K}{\frown} z$ and $z \underset{K}{\frown} f$ with coefficients in $P\otimes_A M$ and $M\otimes_A P$ respectively, as follows.
\\
1. If $p$ and $q-p$ are not both odd, so that $\nu(q-p)=\nu(q)-\nu(p)$, one has
\begin{eqnarray*}
f\underset{K}{\frown} z = (f(x_{\nu(q-p)+1} \ldots  x_{\nu(q)})\otimes_A m)\otimes \, x_1 \ldots x_{\nu(q-p)},\\
z\underset{K}{\frown} f = (-1)^{pq} (m \otimes_A f(x_1 \ldots  x_{\nu(p)}))\otimes \, x_{\nu(p)+1} \ldots x_{\nu(q)}. 
\end{eqnarray*}
2. If $p$ is odd and $q$ is even, so that $\nu(q-p)=\nu(q)-\nu(p)-N+2$, one has 
\begin{eqnarray*}
  f\underset{K}{\frown} z =  -\sum_{0\leq i+j \leq N-2} (x_{\nu(q-p)+i+1}\ldots x_{\nu(q)-\nu(p)-j} f(x_{\nu(q)-\nu(p)-j+1} \ldots x_{\nu(q)-j}) \\
  \otimes_A x_{\nu(q)-j+1}\ldots x_{\nu(q)}m x_1 \ldots x_i)\otimes \, x_{i+1} \ldots  x_{\nu(q-p)+i}.
\end{eqnarray*}
\begin{eqnarray*}
  z\underset{K}{\frown} f =  \sum_{0\leq i+j \leq N-2} (x_{\nu(q)-j+1}\ldots x_{\nu(q)}m x_1 \ldots x_i \otimes_A f(x_{i+1} \ldots x_{i+\nu(p)})\\  x_{i+\nu(p)+1}\ldots x_{\nu(p)+N-j-2}) 
  \otimes \, x_{\nu(p)+N-j-1}\ldots x_{\nu(q)-j}.
\end{eqnarray*}
These formulas hold for $q \geq p$. If $q<p$, one has $f \underset{K}{\frown} z=z \underset{K}{\frown} f=0$.

For our example $A=S(V,N)$, we rewrite these formulas in terms of wedge product as in the previous section. Using the isomorphisms (\ref{isocan}), we see that $f:\bigwedge^{\nu(p)} V \rightarrow P$ and $z=m \otimes x_1 \wedge  \ldots \wedge x_{\nu(q)}$ in $M\otimes \bigwedge^{\nu(q)} V$, where $x_1< \ldots <x_{\nu(q)}$ belong to the basis $(e_1, \ldots , e_n)$ of $V$.

If $\nu(q-p)=\nu(q)-\nu(p)$, we obtain
\begin{eqnarray} \label{leftcap1}
f\underset{K}{\frown} z = (-1)^{(q-p)p} \sum _{\sigma \in Sh(\nu(q-p),\nu(p))} \mathrm{sgn}(\sigma) \,(f(x_{\sigma(\nu(q-p)+1)} \wedge \ldots \wedge x_{\sigma(\nu(q))})\otimes_A m)\nonumber \\
\otimes \, (x_{\sigma(1)}\wedge \ldots \wedge x_{\sigma(\nu(q-p))}).   
\end{eqnarray}

If $\nu(q-p)\neq \nu(q)-\nu(p)$, meaning that $N>2$, $p$ is odd and $q$ is even, we obtain
\begin{eqnarray} \label{leftcap2}
  \lefteqn{f\underset{K}{\frown} z = - \sum_{0\leq i+j \leq N-2} \ \ \sum _{\sigma \in Sh(i,\nu(q-p),\nu(p),j)} \mathrm{sgn}(\sigma) } \nonumber \\
& &   (x_{\sigma (\nu(q-p)+i+1)}\ldots x_{\sigma(\nu(q)-\nu(p)-j)} f(x_{\sigma (\nu(q)-\nu(p)-j+1)} \wedge \ldots \wedge  x_{\sigma(\nu(q)-j)}) \\
& & \otimes_A x_{\sigma(\nu(q)-j+1)}\ldots x_{\sigma(\nu(q))} m x_{\sigma(1)} \ldots x_{\sigma(i)}) \otimes \, (x_{\sigma(i+1)}\wedge \ldots  \wedge x_{\sigma(\nu(q-p)+i)}), \nonumber  
\end{eqnarray}
where $Sh(i,\nu(q-p),\nu(p),j)$ denotes the subset of $\Sigma_n$ formed of the $(\nu(q-p),\nu(p))$-shuffles acting on the index sets $\{i+1, \ldots ,i+\nu(q-p)\}$ and $\{\nu(q)-\nu(p)-j+1, \ldots , \nu(q)-j\}$.

Similarly, if $\nu(q-p)=\nu(q)-\nu(p)$, we obtain 
\begin{eqnarray} \label{rightcap1}
z\underset{K}{\frown} f = (-1)^{qp} \sum _{\sigma \in Sh(\nu(p),\nu(q-p))} \mathrm{sgn}(\sigma) \,(m \otimes_A  f(x_{\sigma(1)} \wedge \ldots \wedge x_{\sigma(\nu(p))}))\nonumber \\
\otimes \, (x_{\sigma(\nu(p)+1)}\wedge \ldots \wedge x_{\sigma(\nu(q))}),   
\end{eqnarray}
while if $\nu(q-p)\neq \nu(q)-\nu(p)$,
\begin{eqnarray} \label{rightcap2}
  \lefteqn{z\underset{K}{\frown} f = \sum_{0\leq i+j \leq N-2} \ \ \sum _{\sigma \in Sh(i,\nu(p),\nu(q-p),j)} \mathrm{sgn}(\sigma) } \nonumber \\
& &   (x_{\sigma (\nu(q)-j+1)}\ldots x_{\sigma(\nu(q))} m x_{\sigma(1)} \ldots x_{\sigma(i)} \otimes_A f(x_{\sigma (i+1)} \wedge \ldots \wedge  x_{\sigma(i+\nu(p))}) \\
& &  x_{\sigma(i+\nu(p)+1)}\ldots x_{\sigma(\nu(p)+N-j-2)}) \otimes \, (x_{\sigma(\nu(p)+N-j-1)}\wedge \ldots  \wedge x_{\sigma(\nu(q)-j)}), \nonumber  
\end{eqnarray}
where $Sh(i,\nu(p),\nu(q-p),j)$ denotes the subset of $\Sigma_n$ formed of the $(\nu(p),\nu(q-p))$-shuffles acting on the index sets $\{i+1, \ldots ,i+\nu(p)\}$ and $\{\nu(p)+N-j-1, \ldots , \nu(q)-j\}$.
\\ 

If $N=2$ and $M=P=A$, the maps $f \underset{K}{\frown} -$ and $- \underset{K}{\frown} f$ are equal up to a sign to the contraction map (inner product) $i_f$ of the Cartan calculus. For the Cartan calculus, the reader may consult~\cite{am:foundmech} Section 2.4, or for a more recent reference,~\cite{lgpv:poisson} Chapter 3.
\\

\textbf{Question 2.} The same as in the previous section, but now for the $N$-Koszul cap products. To resume both questions, is it possible to define an $N$-Cartan calculus on manifolds? If the answer is affirmative, the algebraic structure of the $N$-Koszul calculus~\cite{rb:Ncal} would extend to the so-defined $N$-Cartan calculus on manifolds, and would provide a new prototype of noncommutative differential geometry, different from a Tamarkin-Tsygan calculus~\cite{tt:calculus}.

\vspace{0.5 cm} \textsf{Roland Berger: Univ Lyon, UJM-Saint-\'Etienne, CNRS UMR 5208, Institut Camille Jordan, F-42023, Saint-\'Etienne, France}

\emph{roland.berger@univ-st-etienne.fr}\\


\begin{thebibliography}{99}
  
\bibitem{am:foundmech} R. Abraham, J. E. Marsden, \emph{Foundations of Mechanics}, Second Edition, Advanced Book Program, The Benjamin/Cummings Publishing Company, 1978. 
\bibitem{rb:nonquad} R. Berger, Koszulity for nonquadratic algebras, \emph{J. Algebra} \textbf{239} (2001), 705-734.
\bibitem{rb:Ncal} R. Berger, Koszul calculus for $N$-homogeneous algebras, arXiv:1610.01035.
\bibitem{rbnm:kogo} R. Berger, N. Marconnet, Koszul and Gorenstein properties for homogeneous algebras, 
  \emph{Algebras and Representation Theory} \textbf{1} (2006), 67-97.
\bibitem{lgpv:poisson} C. Laurent-Gengoux, A. Pichereau, P. Vanhaecke, \emph{Poisson Structures}, Grundlehren der mathematischen Wissenschaften \textbf{347}, Springer, 2013.
\bibitem{tt:calculus} D. Tamarkin, B. Tsygan, The ring of differential operators on forms in noncommutative calculus, \emph{Proc. Sympos. Pure Math.} \textbf{73}, Amer. Math. Soc. (2005), 105-131.
\end{thebibliography}
\end{document}